\documentclass[11pt, reqno]{article}
\usepackage{graphicx,indentfirst}
\usepackage{amscd}

\usepackage[all]{xy}
\usepackage{amsmath}
\usepackage{amsfonts,amssymb}
\usepackage{amsthm}
\usepackage{latexsym}
\usepackage[american]{babel}
\usepackage{times}
\usepackage{tikz}
\usepackage{comment}
\usetikzlibrary{matrix,arrows,decorations.pathmorphing}
\usepackage[a4paper,left=2.5cm,right=2.5cm,bottom=3cm,top=3.5cm]{geometry}
\RequirePackage[T1]{fontenc}

\begin{document}
\title{A proof of the Baum-Connes conjecture for real semisimple Lie groups with coefficients on flag
varieties}
\author{Zhaoting Wei\\Kent State University at Geauga\\
14111 Claridon Troy Road, Burton, OH 44021, USA\\zwei3@kent.edu}
\pagestyle{plain}

\newcommand{\End}{\text{End}}
\newcommand{\ad}{\text{ad}}
\newcommand{\Ad}{\text{Ad}}
\newcommand{\Pol}{\text{Pol}}
\newcommand{\Cl}{\text{Cl}}
\newcommand{\B}{\text{B}}
\newcommand{\Image}{\text{Image}}
\newcommand{\K}{\text{K}}
\newcommand{\KK}{\text{KK}}
\newcommand{\SL}{\text{SL}}
\newcommand{\red}{\text{red}}
\newcommand{\pt}{\text{pt}}

\newtheorem{thm}{Theorem}[section]
\newtheorem{lemma}[thm]{Lemma}
\newtheorem{prop}[thm]{Proposition}
\newtheorem{coro}[thm]{Corollary}
\newtheorem{defi}{Definition}[section]
{\theoremstyle{plain} \newtheorem{rmk}{Remark}}
{\theoremstyle{plain} \newtheorem{eg}{Example}}

\maketitle

\begin{abstract}
We consider the equivariant K-theory of a real semisimple Lie group which acts on the (complex) flag variety of its complexification group. We construct an assemble map in the framework of KK-theory and then we prove that it is an isomorphism. The prove relies on a careful study of the orbits of the real group action on the flag variety and then piecing together different orbits. This result is a special case of the Baum-Connes conjecture with coefficients.


\end{abstract}


\section{Introduction}\label{section:introduction}
Let $G$ be a locally compact topological group and $A$ is a $C^*$-algebra equipped with a continuous action of G by $C^*$-algebra automorphisms. Following \cite[Section 4]{block2013weyl}, we define the equivariant K-theory of $A$ to be the K-theory of the \emph{reduced crossed product algebra}:
$$
\K^*_G(A):=\K^*(C^*_r(G,A)).
$$
The equivariant K-theory defined in this way has a useful connection to Baum-Connes conjecture and representation theory.  It is well-known that $C^*_r(G)$ reflects the tempered unitary dual when $G$ is a reductive Lie group, see \cite[Section 4.1]{block2013weyl}.

When $G$ is compact and $\mathcal{X}$ is a compact $G$-topological space, let $C_0(\mathcal{X})$ denote the $C^*$-algebra of compact supported complex value continuous functions on $\mathcal{X}$. Then it is well-known that our $\K^*_G(C_0(\mathcal{X}))$ coincides with the equivariant K-theory of $\mathcal{X}$, see \cite{julg1981k}. If $\mathcal{X}$ is itself compact, we can also denote $C_0(\mathcal{X})$ by $C(\mathcal{X})$.

\begin{rmk}
Be aware that this is not the same as Kasparov's definition of equivariant K-theory in \cite{kasparov1988equivariantkk}.
\end{rmk}

Back to general $G$. For any  $C^*$-algebra $A$ and $B$ with continuous $G$-action we have the equivariant KK-theory group $\KK^G_*(A,B)$ as in \cite[Definition 2.3]{kasparov1988equivariantkk}.

The Baum-Connes conjecture can be formulated as follows. Let $G$ be almost connected and $U$ be its maximal compact subgroup and $\mathcal{S}:=G/U$ be the quotient space.  We have the \emph{assemble map} \cite{higson1998baum}
$$
\mu_{G}: \KK^G_*(C_0(\mathcal{S}),\mathbb{C})\rightarrow \K_G^*(\mathbb{C}).
$$
The \emph{Baum-Connes conjecture} claims that the assemble map $\mu_{G}$ is an isomorphism. In 2003, J. Chabert, S. Echterhoff, R. Nest \cite{chabert2003connes} proved this conjecture for almost connected Lie groups and for linear $p$-adic groups. 

Recall that a topological group $G$ is called \emph{almost connected} if $G/G_0$ is compact, where $G_0$ denotes the identity component of $G$. For a Lie group $G$, almost connected simply means that $G$ has finitely many connected components.

\begin{rmk}
Baum-Connes conjecture is still open for certain discrete groups, for example $G=\SL(3,\mathbb{Z})$. Notice that in this case $G$ is not almost connected and the space $\mathcal{S}=G/U$ should be replaced by the \emph{universal proper $G$-space} $\underline{E}G$. See \cite{baum1994classifying} for more details. For an introduction of Baum-Connes conjecture for discrete groups we refer to \cite{valette2002introduction}.
\end{rmk}

Moreover, for a  $C^*$-algebra $A$ with continuous $G$-action, we have the Baum-Connes conjecture conjecture with coefficients in $A$, which claims that the map
$$
\mu_{G,A}: \KK^G_*(C_0(\mathcal{S}),A)\rightarrow \K_G^*(A)
$$
is an isomorphism. There are counter examples for some certain $G$ and $A$ as in  \cite{higson2002counterexamples}. For general $G$ and $A$ Baum-Connes conjecture with coefficients is still  open.

In this paper we focus on the case that $G$ is a real semisimple Lie group and $A$ is the  $C^*$-algebra of continuous functions on the complexified flag variety of $G$. In more details let $G_{\mathbb{C}}$ be the complexification of $G$, We have the flag variety $\mathcal{B}$ of $G_{\mathbb{C}}$. The group $G_{\mathbb{C}}$ (hence $G$ and $U$) acts on $\mathcal{B}$, so we also have the assemble map
\begin{equation}\label{defiofmap}
\mu_{G,\mathcal{B}}: \KK^G_*(C_0(\mathcal{S}),C(\mathcal{B}))\rightarrow \K_G^*(C(\mathcal{B})).
\end{equation}

The main result of this paper is the following theorem:
\begin{thm}\label{Baum Connes flag}
For any real semisimple Lie group $G$, the assemble map
$$
\mu_{G,\mathcal{B}}: \KK^G_*(C_0(\mathcal{S}),C(\mathcal{B}))\rightarrow \K_G^*(C(\mathcal{B}))
$$
is an isomorphism.
\end{thm}

\begin{rmk}
The significance of $\K_G^*(C(\mathcal{B}))$ has been discussed in \cite[Section 4.4]{block2013weyl}.
\end{rmk}

The proof of Theorem \ref{Baum Connes flag}  in this paper relies on a careful study of the orbits of the real group action on the flag variety: We first proof the isomorphism on one single orbit of the $G$-action by reducing to solvable subgroups, and then we piece together  assemble maps on different orbits. The proof does not require the hard techniques in functional analysis and representation theory so it can be considered as an \emph{geometric proof}.

This paper is organized as follows: In Section \ref{section: flag varieties} and \ref{section: KK-theory and the assemble map} we construct the assemble map. In Section \ref{section:The assemble Map on a Single G-Orbit of F} we study the assemble map on one single $G$-orbit of the flag variety. In Section \ref{section:The G-orbits on the Flag Variety} we study the $G$-orbits on $\mathcal{B}$ and in Section \ref{section:Baum-Connes Conjecture on Flag Varieties} we prove the Theorem \ref{Baum Connes flag}. In Section \ref{section:An Example} we give an example to illustrate the idea of the construction.

This work is inspired by the study of equivariant K-theory in \cite{block2013weyl} and  \emph{Matsuki correspondence} in \cite{mirkovic1992matsuki}.

\subsection*{Acknowledgement}
The author would like to thank Jonathan Block and Nigel Higson for very helpful discussions and comments. The author also wants to thank Weibo Fu, Yuhang Liu, and Jun Su for kindly answering questions related to this work.

\section{Real Semisimple Lie Groups and Flag Varieties}\label{section: flag varieties}
We will use the following notations in this paper. Let $G$ be a connected linear real semisimple Lie group, $U$ be the identity component of a maximal compact subgroup of $G$. In the sequel we fix such a $U$ and call it \emph{the} maximal compact subgroup of $G$. We denote the space $G/U$ by $\mathcal{S}$.

Let $G_{\mathbb{C}}$ be the complexification of $G$, $B_{\mathbb{C}}$ be the Borel subgroup of $G_{\mathbb{C}}$ and $\mathcal{B}=G_{\mathbb{C}}/B_{\mathbb{C}}$ be the flag variety.

Obviously $G$ acts on the flag variety $\mathcal{B}$. Unlike $G_{\mathbb{C}}$, the $G$-action is not transitive, see  Section \ref{section:The G-orbits on the Flag Variety} below or \cite{mirkovic1992matsuki}.

\begin{eg}\label{eg: sl2}
If $G=\SL(2,\mathbb{R})$ then $U=SO(2)$ and $\mathcal{S}=G/U=\mathbb{H}$ the upper half plane.

On the other hand $G_{\mathbb{C}}=\SL(2,\mathbb{C})$. Hence
$$
B_{\mathbb{C}}=\left\{\left.\begin{pmatrix}a &b \\ 0 & a^{-1}\end{pmatrix}\right|a\in \mathbb{C}^*, B\in\mathbb{C}\right\}
$$
and
$$
\mathcal{B}=G_{\mathbb{C}}/B_{\mathbb{C}}=\mathbb{C}P^1\cong S^2.
$$
$G_{\mathbb{C}}$ (hence $G$) acts on $\mathcal{B}=\mathbb{C}P^1$ by fractional linear transform. In projective coordinates we have
$$
\begin{pmatrix}a &b \\ c & d\end{pmatrix}\cdot \begin{pmatrix}u \\ v\end{pmatrix}:=\begin{pmatrix}au+bv \\ cu+dv\end{pmatrix}.
$$
If we set $z=u/v$, then
$$
\begin{pmatrix}a &b \\ c & d\end{pmatrix}\cdot z:= \frac{az+b}{cz+d}.
$$
We will study $G$-orbits of $\mathcal{B}$ in more details in Section \ref{section:An Example}.
\end{eg}

\section{KK-theory and the Assemble Map}\label{section: KK-theory and the assemble map}

In this section we quickly review KK-theory and construct the assemble map
\begin{equation}
\mu_{G,\mathcal{T}}: \KK^G_*(\mathcal{S},\mathcal{T})\rightarrow \K_G^*(\mathcal{T})
\end{equation}
for any $G$-space $\mathcal{T}$. We work in the framework of Kasparov as in \cite{kasparov1988equivariantkk}. 

In this paper we use KK-theory as a black box and most results in this section are given without proof.

\subsection{A quick review of equivariant K-theory}
Let  $G$ be a locally compact group. We call a $C^*$-algebra with continuous $G$ action a \emph{$G$-$C^*$-algebra}. For a $G$-$C^*$-algebra $A$,  we define the \emph{reduced cross product $C^*$-algebra} $C^*_r(G,A)$ as the completion of the twisted convolution algebra of compactly supported and continuous functions from $G$ into $A$. The convolution product is
$$
f_1\star f_2(g)=\int_Gf_1(h)\alpha_h(f_2(h^{-1}g))dh
$$
where $\alpha$ denotes the action of $G$ on $A$. If $A$ is represented faithfully and
isometrically on a Hilbert space $\mathcal{H}$, then the completion is under the operator norm on $L^2(G,\mathcal{H})$. In particular, $C^*_r(G,\mathbb{C})= C^*_r(G)$. See \cite[Chapter 7]{eilers2018c} for details.

 Let $\mathcal{T}$ be a topological space with continuous $G$-action. Let $C_0(\mathcal{T})$ be the space of continuous functions on $\mathcal{T}$ which vanishes at infinity. If $\mathcal{T}$ is compact, then $C_0(\mathcal{T})=C(\mathcal{T})$ is the space of all continuous functions on $\mathcal{X}$. We define
    \begin{equation}
    \K^*_G(\mathcal{T}):=\K^*(C^*_r(G,C_0(\mathcal{T}))).
    \end{equation}

\begin{itemize}
 \item $\K^*_G(\pt)$ reflects the tempered unitary dual when $G$ is a reductive Lie group, see \cite{baum1994classifying}.

\item When $G$ is compact and $\mathcal{X}$ is a compact $G$-topological space, let $C(\mathcal{X})$ denote the $C^*$-algebra of complex value continuous functions on $\mathcal{X}$. Then our $\K^*_G(\mathcal{X})$ coincides with the equivariant K-theory of $\mathcal{X}$, see \cite{julg1981k}.
\end{itemize}

\subsection{A quick review of equivariant KK-theory}
For two $G$-$C^*$-algebras $A$ and $B$, Kasparov introduced the equivariant KK-theory $\KK^G_i(A,B)$ for $i=0, 1$ in \cite{kasparov1988equivariantkk}. In this paper we do not go to details of the construction but we list some properties of KK-theory here.

\begin{prop}
$\KK^G_i(A,B)$ is an abelian group for $i=0, 1$, and it is contravariant for $A$ and covariant for $B$. 
\end{prop}

Recall that a $C^*$-algebra  is called \emph{$\sigma$-unital} if it possesses a countable approximate unit.

\begin{thm}\label{thm: Kasparov product}[\cite[Theorem 2.11, Definition 2.12, and Theorem 2.14]{kasparov1988equivariantkk}]
If $A$ is separable, then we have the r \emph{Kasparov product}
$$
\KK^G_i(A,B_1)\otimes \KK^G_j(B_1,B_2)\to KK^G_{i+j}(A,B_2).
$$
More generally if $A_1$ and $A_2$ are separable, then we have the Kasparov product
\begin{equation}\label{eq: Kasparov product}
\KK^G_i(A_1,B_1\hat{\otimes}D)\otimes \KK^G_j(D\hat{\otimes}A_2,B_2)\to KK^G_{i+j}(A_1\hat{\otimes}A_2,B_1\hat{\otimes}B_2).
\end{equation}
denoted by $x_1\otimes_D x_2$. Moreover, the Kasparov product has the following properties
\begin{enumerate}
\item It is bilinear;
\item It is contravariant in $A_1$ and $A_2$ and covariant in $B_1$ and $B_2$;
\item It is functorial in $D$;
\item It is associative;
\item For any $\sigma$-unital $G$-$C^*$-algebra $A$, there exists a two side multiplicative unit $1_A\in \KK^G_0(A,A)$.
\end{enumerate}
\end{thm}

\begin{prop}\label{prop: sigmaD}
For a $\sigma$-unital $G$-$C^*$-algebra $D$, we have a homomorphism
\begin{equation}
\sigma_D:\KK^G_i(A,B) \to \KK^G_i(A\hat{\otimes}D,B\hat{\otimes}D),
\end{equation}\label{eq: sigma D}
 where $\hat{\otimes}$ denotes the tensor product completed under the minimal norm. The map $\sigma_D$ is compatible with the Kasparov product in the sense that if $A_1$, $A_2$, and $D_1$ are separable, then
$$
\sigma_{D_1}(x_1\otimes_D x_2)=\sigma_{D_1}(x_1)\otimes_{D\hat{\otimes}D_1}\sigma_{D_1}(x_2) 
$$
for $x_1\in \KK^G_*(A_1,B_1\hat{\otimes}D)$ and $x_2\in \KK^G_*(D\hat{\otimes}A_2,B_2)$.
\end{prop}

\begin{prop}\label{prop: restriction of KK wrt groups}
Let $f: G_1\rightarrow G_2$ be a homomorphism between groups, we have the natural \emph{restriction homomorphism}
$$
r^{G_2,G_1}: \KK^{G_2}_*(A, B)\longrightarrow \KK^{G_1}_*(A, B)
$$
which is compatible with the Kasparov product.
\end{prop}

\begin{prop} \label{fromequiKKtocrossedproductKK}[\cite[Theorem 3.11]{kasparov1988equivariantkk}]
There is a natural homomorphism
 $$
 j^G_r: \KK^G_*(A,B)\longrightarrow \KK_*(C^*_r(G,A),C^*_r(G,B))
 $$

which is compatible with the Kasparov product. Here $\KK_*(-,-)$ denotes the ordinary (non-equivariant) KK-theory. Moreover, for $1_A\in \KK^G_0(A,A)$ we have
$$
 j^G_r(1_A)=1_{C^*_r(G,A)}\in \KK_0(C^*_r(G,A),C^*_r(G,A)).
$$
\end{prop}

As before, if $A=C_0(\mathcal{X})$ and $B=C_0(\mathcal{Y})$ for topological spaces $\mathcal{X}$ and $\mathcal{Y}$, then we denote $\KK^G_i(C_0(\mathcal{X}), C_0(\mathcal{Y}))$ simply by $\KK^G_i(\mathcal{X}, \mathcal{Y})$.

We have the \emph{Poincare duality} isomorphism in KK-theory.
\begin{thm}\label{Poincare duality in KK theory}[\cite[Theorem 4.10]{kasparov1988equivariantkk}, see also \cite[Section 4.3]{block2013weyl}]
For a $G$-manifold $\mathcal{X}$, let $C_{\tau}(\mathcal{X})$ denote the algebra of continuous sections of the Clifford bundle over $\mathcal{X}$ vanishing at infinity. Then we have the following isomorphism
\begin{equation}
\KK^G_*(\mathcal{X},\mathcal{T})\cong \K_G^*(C_0(\mathcal{T})\otimes C_{\tau}(\mathcal{X})).
\end{equation}
\end{thm}

\subsection{The Dirac Element}
For $G$, $\mathcal{X}$ and $C_{\tau}(\mathcal{X})$ as in Theorem \ref{Poincare duality in KK theory}, Kasparov defined the \emph{Dirac element} \cite[Section 4.2]{kasparov1988equivariantkk}:
\begin{equation}
\mathfrak{d}_{G,\mathcal{X}}\in \KK^G_0(C_{\tau}(\mathcal{X}),\mathbb{C})
\end{equation}

\begin{rmk}\label{Diracnotspin}
In the definition of $\mathfrak{d}_{G,\mathcal{X}}$ we do not require that $\mathcal{X}$ is spin. We will discuss the spin case in Section \ref{subsection: spin case} below.
\end{rmk}

Now we want to find the relation between equivariant KK-theory and the K-theory of crossed-product algebras.

First remember that for any $G$-space $\mathcal{T}$ we have the map as in \eqref{eq: sigma D}
$$
\sigma_{\mathcal{T}}: \KK^G_*(A,B)\longrightarrow \KK^G_*(A\otimes C_0(\mathcal{T}),B\otimes C_0(\mathcal{T})).
$$
Apply $\sigma_{\mathcal{T}}$ to $\mathfrak{d}_{G,\mathcal{S}}\in \KK^G_0(C_{\tau}(\mathcal{S}),\mathbb{C})$ we get
$$
\sigma_{\mathcal{T}}(\mathfrak{d}_{G, \mathcal{S}})\in \KK^G_0(C_{\tau}(\mathcal{S})\otimes C_0(\mathcal{T}),C_0(\mathcal{T})).
$$
Then apply the map $ j^G_r$ in Proposition \ref{fromequiKKtocrossedproductKK} to $\sigma_{\mathcal{T}}(\mathfrak{d}_{G, \mathcal{S}})$ we get
$$
 j^G_r(\sigma_{\mathcal{T}}(\mathfrak{d}_{G, \mathcal{S}}))\in \KK_0(C^*_r(G,C_{\tau}(\mathcal{S})\otimes C_0(\mathcal{T})),C^*_r(G,C_0(\mathcal{T}))).
$$
We denote $ j^G_r(\sigma_{\mathcal{T}}(\mathfrak{d}_{G,\mathcal{S}}))$ by $\mathfrak{D}_{G,\mathcal{S}}$ or simply by $\mathfrak{D}$ if $G$ is clear.

\begin{defi}[The assemble map]\label{assemble map}
 Let $\mathcal{S}=G/U$, for any $\mathcal{T}$, the Poincare duality and the Kasparov product with $\mathfrak{D}$ give us the desired map
\begin{equation}
\cdot \otimes \mathfrak{D}: \KK^G_*(\mathcal{S},\mathcal{T})\cong \K_G^*(C_0(\mathcal{T})\otimes C_{\tau}(\mathcal{S}))\rightarrow \K_G^*(\mathcal{T}).
\end{equation}
\end{defi}

\begin{rmk}\label{assemblenotspin}
As pointed out in Remark \ref{Diracnotspin}, we do not require $S$ to be spin to define the assemble map.
\end{rmk}

\subsection{The Spin Case}\label{subsection: spin case}

Let us study the assemble map in the spin case to get more intuition.

When $\mathcal{S}$ is spin and   even dimensional, it is well known that $C_{\tau}(\mathcal{S})$ is strongly Morita equivalent to $C_0(\mathcal{S})$. Hence the Poincare duality gives us
\begin{equation}
\KK^G_*(\mathcal{S},\mathcal{T})\cong \K_G^*(\mathcal{T} \times \mathcal{S}).
\end{equation}

In this case, the Dirac element $\mathfrak{d}_{G,\mathcal{S}}$ is exactly the index map of the Dirac operator in the $\mathcal{S}$ direction (\cite{atiyah1968bott}) and this justified the name "Dirac element". In this case the assemble map is given by the index map
\begin{equation}
\mathfrak{D}: \K_G^*(\mathcal{T} \times \mathcal{S})\rightarrow \K_G^*(\mathcal{T})
\end{equation}

We can look at $\K_G^*(\mathcal{T} \times \mathcal{S})$ from another viewpoint. Remember that $\mathcal{S}=G/U$. We have the following obvious result
\begin{lemma}\label{changefiberproducttocartesianproduct}
Let $\mathcal{T}$ be a $G$-space in the above setting and $H$ be a subgroup of $G$. Then $G \times_H \mathcal{T}$ is $G$-isomorphic to $G/H \times \mathcal{T}$, where $G$ acts on $G/H \times \mathcal{T}$ by the diagonal action. Hence
$$
\K^*_G(G/H \times \mathcal{T})\cong \K^*_G(G \times_H \mathcal{T}).
$$
 \end{lemma}
\begin{proof} The  map
\begin{align*}
G\times_H \mathcal{T}\rightarrow & G/H \times \mathcal{T}\\
(g, t)\mapsto & (g, gt)
\end{align*}
gives the $G$-isomorphism
\end{proof}

We also have the following isomorphism, see \cite{rieffel1982applications}
\begin{lemma}[The induction map]\label{Ktheoryinductionmap}
Consider a group $G$ and  a closed subgroup $H\subset G$. For  an $H$-space $\mathcal{T}$, there is an induction map
$$
\K^*_H(\mathcal{T})\rightarrow \K^*_G(G \times_H \mathcal{T})
$$
which is a natural isomorphism. Here the $G$-action on $G \times_H \mathcal{T}$ is the left multiplication on the first component.
\end{lemma}
\begin{proof}Just notice that $C^*_r(H, C_0(\mathcal{T}))$ and $C^*_r(G, C_0(G \times_H \mathcal{T}))$ are strongly Morita equivalent.\end{proof}

\begin{rmk}
Although simple, the idea of Lemma \ref{changefiberproducttocartesianproduct} and \ref{Ktheoryinductionmap} will appear later in Lemma \ref{interchangesubgroups}.
\end{rmk}

\begin{coro}\label{coro: reduce to max compact subgroup}
Let $G$ be an almost connected Lie group and $U$ be its maximal compact subgroup. If $\mathcal{S}=G/U$ is spin and even dimensional, then for any $G$-space $\mathcal{T}$ we have a natural isomorphism
\begin{equation}
\KK^G_*(\mathcal{S},\mathcal{T})\cong \K_G^*(\mathcal{T} \times \mathcal{S}) \cong K_U^*(\mathcal{T}).
\end{equation}
\end{coro}
 
 According to Corollary \ref{coro: reduce to max compact subgroup}, the assemble map in Definition \ref{assemble map} has the following form
\begin{equation}
\mathfrak{D}: \K_U^*( \mathcal{T})\rightarrow \K_G^*( \mathcal{T}).
\end{equation}
The \emph{Connes-Kasparov conjecture}, which is a special case of the Baum-Connes conjecture, claims that the above map is an isomorphism.

\begin{rmk}
The original Connes-Kasparov conjecture does not require $G/U$ to be spin but it is stated in a slightly different way, see \cite{penington1983dirac}.
\end{rmk}

\subsection{The Dual Dirac Element}\label{subsection: dual Dirac}
We are looking for an inverse element of $\mathfrak{d}_{G,\mathcal{X}}\in \KK^G_0(C_{\tau}(\mathcal{X}),\mathbb{C})$. For this purpose Kasparov  introduced the concept of $G$-\emph{special manifold} in \cite[Section 5.1]{kasparov1988equivariantkk}. 
\begin{defi}\label{defi: G-special manifolds}
A $G$- manifold $\mathcal{X}$ is called $G$-special if there exists an element $
 \eta_{G,\mathcal{X}}\in  \KK^G_0(\mathbb{C},C_{\tau}(\mathcal{X}))
 $ called the \emph{dual Dirac element}, such that 
 $$
 \mathfrak{d}_{G,\mathcal{X}}\otimes_\mathbb{C}\eta_{G,\mathcal{X}}=1_{C_{\tau}(\mathcal{X})}
 $$
under the Kasparov product $
\KK^G_0(C_{\tau}(\mathcal{X}),\mathbb{C})\otimes \KK^G_0(\mathbb{C},C_{\tau}(\mathcal{X}))\rightarrow \KK^G_0(C_{\tau}(\mathcal{X}),C_{\tau}(\mathcal{X}))
 $.
\end{defi}

\begin{rmk}
It is clear that the element $\eta_{G,\mathcal{X}}$ is unique if exists.
\end{rmk}

\begin{defi}\label{defi: gamma X}
 We consider $\eta_{G,\mathcal{X}} \otimes_{C_{\tau}(\mathcal{X})} \mathfrak{d}_{G,\mathcal{X}}\in \KK^G_0(\mathbb{C},\mathbb{C})$ 
as the Kasparov product in the other way and we denote it by $\gamma_{G,\mathcal{X}}$.
\end{defi}

\begin{rmk}
Even for a $G$-special manifold $\mathcal{X}$, the element $
 \gamma_{G,\mathcal{X}} \in \KK^G_0(\mathbb{C},\mathbb{C})
 $
 need not to be $1_{\mathbb{C}}$. If $\gamma_{G,\mathcal{X}}=1_{\mathbb{C}}$, then $\mathfrak{d}_{G,\mathcal{X}}$ and $\eta_{G,\mathcal{X}}$ are inverse to each other under the Kasparov product.
 \end{rmk}

In \cite{kasparov1988equivariantkk} Kasparov gave several examples of $G$-special manifolds, in particular he gave the following result.
\begin{lemma}\label{lemma: G/U is G-special}[\cite[Theorem 5.7]{kasparov1988equivariantkk}]
Let $G$ be an almost connected group and $U$ the maximal compact subgroup, then  the homogeneous space $\mathcal{S}=G/U$ is a $G$-special manifold. Moreover the element $\gamma_{G,\mathcal{S}}$ is independent of the choice of $U$.
\end{lemma}

\subsection{When $\gamma_{G,\mathcal{X}}=1$?}
Kasparov proved that $\gamma_{G,\mathcal{X}}=1$ in some special cases, which is sufficient for our purpose. First we recall the concept of amenable group.

\begin{defi}\label{defi: amenable group}
A group $G$ is called \emph{amenable} if there exists a left invariant mean $\mu$ on the space $L^{\infty}(G)$. Here we call $\mu: L^{\infty}(G)\to \mathbb{R}$ a mean if it is a non-negative linear functional such that $\mu(1)=1$. A mean $\mu$ is called left invariant if
$\mu(f)=\mu(l_g(f))$ for any $f\in L^{\infty}(G)$ and $g\in G$, where $l_gf(g^{\prime}):=f(gg^{\prime})$.
\end{defi}

\begin{eg}
\begin{itemize}
\item An abelian group is amenable;
\item A solvable group is amenable;
\item A compact group is amenable;
\item Every (closed) subgroup of an amenable group is amenable;
\item A non-compact real (or complex) semisimple Lie group is never amenable.
\end{itemize}
\end{eg}

\begin{lemma}\label{restrictionofDiracelements}
If $\mathcal{X}$ is a $G_2$-special manifold and $f: G_1\rightarrow G_2$ be a homomorphism of groups. Then $\mathcal{X}$ is also $G_1$-special and we have
$r^{G_2,G_1}(\mathfrak{d}_{G_2,\mathcal{X}})=\mathfrak{d}_{G_1,\mathcal{X}}$, 
$r^{G_2,G_1}(\eta_{G_2,\mathcal{X}})=\eta_{G_1,\mathcal{X}}$, and
$$
r^{G_2,G_1}(\gamma_{G_2,\mathcal{X}})=\gamma_{G_1,\mathcal{X}},
$$
where $r^{G_2,G_1}$ is the restriction homomorphism as in Proposition \ref{prop: restriction of KK wrt groups}.
\end{lemma}
\begin{proof}
It is clear since $r^{G_2,G_1}$ is compatible with the Kasparov product.
\end{proof}

Kasparov proved the following result for amenable groups.

\begin{thm}\label{Amemnablegammaelementforkernel}[\cite[Theorem 5.9]{kasparov1988equivariantkk}]
Let $f: G_1\rightarrow G_2$ be a homomorphism between almost connected groups with the kernel $\ker f$ amenable and the image closed. Let $U_i$ be the maximal subgroup of $G_i$ and $\mathcal{S}_i=G_i/U_i$ for $i=1,2$. Without loss of generality we assume $f(U_1)\subset U_2$. Then the restriction homomorphism gives us
\begin{equation}
r^{G_2,G_1}(\gamma_{G_2,\mathcal{S}_2})=\gamma_{G_1,\mathcal{S}_1}.
\end{equation}
\end{thm}

\begin{coro}\label{coro: restriction of gamma}
For an almost connected group $G$, let $H<G$ be a closed subgroup. Without loss of generality we choose the maximal compact subgroup $U$ of $G$ such that $U\cap H$ is the maximal subgroup of $H$. Let $\mathcal{S}=G/U$ and $\mathcal{S}_H=H/U\cap H$.Then we have
\begin{equation}
r^{G,H}(\gamma_{G,\mathcal{S}})=\gamma_{H,\mathcal{S}_H}.
\end{equation}
\end{coro}

\begin{coro}\label{Amemnablegammaelement}
Let $P$ be an amenable almost connected group and $L$ be the maximal compact subgroup of $P$.  Then
$\gamma_{P,P/L}=1$ hence  $\mathfrak{d}_{P, P/L}$ and $\eta_{P, P/L}$ are inverse to each other in the KK-groups.
\end{coro}

Now we can immediately get an isomorphic result in the almost connected amenable case. The following result is implicitly given in \cite[Section 5.10]{kasparov1988equivariantkk}.
\begin{coro}\label{Diracforamenablegroups}
If $P$ is an almost connected amenable group, $L$ is the maximal compact subgroup of $P$. Then for any $P$-space $\mathcal{T}$, the assemble map
\begin{equation}
\mu_{P,\mathcal{T}}: \KK^P_*(P/L,\mathcal{T})\rightarrow \K_P^*(\mathcal{T})
\end{equation}
is an isomorphism.
\end{coro}
\begin{proof}  By Definition \ref{assemble map}, the assembly map is given by right multiplication with  the element $\mathfrak{D}=j^P_r(\sigma_{\mathcal{T}}(\mathfrak{d}_{P,P/L}))$.
 Corollary \ref{Amemnablegammaelement} tells us that $\mathfrak{d}_{P,P/L}$ is invertible,  and by Proposition \ref{prop: sigmaD} and Proposition \ref{fromequiKKtocrossedproductKK}, both $j^P_r$ and $\sigma_{\mathcal{T}}$ are compatible with the Kasparov product. So $\mu_{P,\mathcal{T}}$ is an isomorphism.
\end{proof}

\begin{rmk}
The results of Corollary \ref{Amemnablegammaelement} and Corollary \ref{Diracforamenablegroups} do not hold for a general  group $G$. So we cannot apply the Dirac-dual Dirac method to prove Baum-Connes conjecture with coefficients in general. Nevertheless in this paper we consider the case that $G$ is a real semisimple Lie group and $\mathcal{T}=\mathcal{B}$ is the flag variety of $G_{\mathbb{C}}$. Although in this case $G$ is not amenable, we can use geometric trick to reduce to the amenable case.
\end{rmk}

\section{The Assemble Map on a Single $G$-Orbit of the Flag Variety}\label{section:The assemble Map on a Single G-Orbit of F}
As in Section \ref{section: flag varieties} let $G$ be a connected real semisimple Lie group and $U$ be its maximal compact subgroup. Let $\mathcal{S}=G/U$ and  $\mathcal{B}$ be the flag variety of $G_{\mathbb{C}}$. It is clear that the $G$-action on $\mathcal{B}$ is not transitive and let us denote $\mathcal{O}_{\alpha}^+$ to be one of the $G$-orbits.

\begin{rmk}
This notation will be justified in Section \ref{section:The G-orbits on the Flag Variety}.
\end{rmk}

The following proposition is the main result of this section.

\begin{prop}\label{prop: AssembleforGUrealssorbit}
Let $G$ be a connected real semisimple Lie group and $U$ be its maximal compact subgroup. Let $\mathcal{S}=G/U$ and  $\mathcal{B}$ be the flag variety of $G_{\mathbb{C}}$. Let  $\mathcal{O}_{\alpha}^+$ be a $G$-orbit in  $\mathcal{B}$. Then the assemble map
\begin{equation}
\mu_{G,\mathcal{O}_{\alpha}^+}: \KK^G_*(\mathcal{S},\mathcal{O}_{\alpha}^+)\rightarrow \K_G^*(\mathcal{O}_{\alpha}^+).
\end{equation}
is an isomorphism.
\end{prop}

\begin{rmk}
Proposition \ref{prop: AssembleforGUrealssorbit}, which focuses on a single $G$-orbit, is the building block of Theorem \ref{Baum Connes flag}. We will piece together different orbits in Section \ref{section:Baum-Connes Conjecture on Flag Varieties}.
\end{rmk}

The proof of Proposition \ref{prop: AssembleforGUrealssorbit} consists of several steps. First we prove the following lemma:
\begin{lemma}[Interchange subgroups]\label{interchangesubgroups}
Let $G$, $\mathcal{S}$, and $\mathcal{O}_{\alpha}^+$ be as in Proposition \ref{prop: AssembleforGUrealssorbit}. Let $H$ be the isotropy group of $G$ at a point point $x \in \mathcal{O}_{\alpha}^+$. Then there is an isomorphism:
$$
\KK^G_*(\mathcal{S}, \mathcal{O}_{\alpha}^+)\overset{\sim}\longrightarrow \KK^H_*(\mathcal{S},\pt).
$$
\end{lemma}
\begin{proof}First by Poincare duality
$$
\KK^G_*(\mathcal{S}, \mathcal{O}_{\alpha}^+)\cong\K_G^*(C_0(\mathcal{O}_{\alpha}^+)\otimes C_{\tau}(\mathcal{S})).
$$

Then notice that $\mathcal{O}_{\alpha}^+$ can by identified with $G/H$. By a strong Morita equivalence argument similar to Lemma \ref{Ktheoryinductionmap} we have
$$
\K_G^*(C_0(\mathcal{O}_{\alpha}^+)\otimes C_{\tau}(\mathcal{S}))\cong \K_H^*(C_{\tau}(\mathcal{S})).
$$

Finally by Poincare duality again we have
$$
\K_H^*(C_{\tau}(\mathcal{S}))\cong\KK^H_*(\mathcal{S},\pt).
$$ We finish the proof. \end{proof}

Next we proof the following result.

\begin{prop}\label{reducetoH}
Let $G$, $\mathcal{S}$, $\mathcal{O}_{\alpha}^+$, and $H$ be as in Proposition \ref{prop: AssembleforGUrealssorbit} and Lemma \ref {interchangesubgroups}. We have the following commutative diagram:
\begin{equation}
\begin{CD}
\KK^G_*(\mathcal{S}, \mathcal{O}_{\alpha}^+)@>\sim>>\KK^H_*(\mathcal{S},\pt)\\
@VV\mu_{G,\mathcal{O}_{\alpha}^+} V                              @VV\mu_{H,\pt} V \\
\K_G^*(\mathcal{O}_{\alpha}^+)           @>\sim>>\K_H^*(pt)
\end{CD}
\end{equation}
where the vertical maps are the assemble maps and the horizontal isomorphisms are given in Lemma \ref{Ktheoryinductionmap} and Proposition \ref{reducetoH}.
\end{prop}

\begin{proof} To prove the proposition we need to have a closer look at the maps. First we look at the right vertical map.
At the beginning we have the Dirac element
$$
\mathfrak{d}_{G,\mathcal{S}}\in \KK^G_0(C_{\tau}(\mathcal{S}),\mathbb{C})
$$
and apply the restriction homomorphism $r^{G,H}$ in Proposition \ref{prop: restriction of KK wrt groups} we get
$$
r^{G,H}(\mathfrak{d}_{G,\mathcal{S}})\in \KK^H_0(C_{\tau}(\mathcal{S}),\mathbb{C})
$$
which by definition equals to $\mathfrak{d}_{H,\mathcal{S}}\in \KK^H_0(C_{\tau}(\mathcal{S}),\mathbb{C})$, the Dirac element of $H$.

Then we apply the map 
$$
j^H_r: \KK^H_*(C_{\tau}(\mathcal{S}),\mathbb{C})\longrightarrow \KK_*(C^*_r(H,C_{\tau}(\mathcal{S})),C^*_r(H))
$$
in Proposition \ref{fromequiKKtocrossedproductKK} and get
$$
j^H_r(r^{G,H}\mathfrak{d}_{G,\mathcal{S}})\in \KK_0(C^*_r(H,C_{\tau}(\mathcal{S})),C^*_r(H))
$$
and we denote it by $\mathfrak{D}_H$. Right multiplication of $\mathfrak{D}_H$ gives the vertical map on the right in the diagram
$$
\KK^H_*(\mathcal{S},\pt)\overset{\mu_{H,\pt}}{\longrightarrow}\K_H^*(\pt).
$$

On the other hand we have the map in Proposition \ref{prop: sigmaD}
$$
\sigma_{\mathcal{O}_{\alpha}^+}: \KK^G_*(C_{\tau}(\mathcal{S}),\mathbb{C})\longrightarrow \KK^G_*(C_{\tau}(\mathcal{S})\otimes C_0(\mathcal{O}_{\alpha}^+),C_0(\mathcal{O}_{\alpha}^+))
$$
so we get
$$
\sigma_{\mathcal{O}_{\alpha}^+}(\mathfrak{d}_{G,\mathcal{S}}) \in\KK^G_0(C_{\tau}(\mathcal{S})\otimes C_0(\mathcal{O}_{\alpha}^+),C_0(\mathcal{O}_{\alpha}^+))
$$
then via $j^G_r$ we get
$$
j^G_r(\sigma_{\mathcal{O}_{\alpha}^+}(\mathfrak{d}_{G,\mathcal{S}}))\in \KK(C^*_r(G,C_{\tau}(\mathcal{S})\otimes C_0(\mathcal{O}_{\alpha}^+)),C^*_r(G,C_0(\mathcal{O}_{\alpha}^+)))
$$
which we denote by $\mathfrak{D}_{G,\mathcal{O}_{\alpha}^+}$. Right multiplication of $\mathfrak{D}_{G,\mathcal{O}_{\alpha}^+}$ gives the left vertical map
$$
\KK^G_*(\mathcal{S},\mathcal{O}_{\alpha}^+)\overset{\mu_{G,\mathcal{O}_{\alpha}^+}}{\longrightarrow}\K_G^*(\mathcal{O}_{\alpha}^+).
$$

The horizontal maps in the diagram are given by strong Morita equivalence. We also notice that under strong Morita equivalence, $\mathfrak{D}_{G,\mathcal{O}_{\alpha}^+}\cong \mathfrak{D}_H$, so the diagram commutes.
\end{proof}

\begin{lemma}\label{lemma: H is amenable}
Let $G$ be a connected real semisimple Lie group and $\mathcal{B}$ be the flag variety of $G_{\mathbb{C}}$. Let $H$ be the isotropy group of $G$ at any point point $x \in \mathcal{B}$ is amenable and almost connected.
\end{lemma}
\begin{proof}
It is clear since $H$ is a closed subgroup of a Borel subgroup of $G_{\mathbb{C}}$.
\end{proof}

According to Proposition \ref{reducetoH}, in order to prove the claim of Proposition \ref{prop: AssembleforGUrealssorbit}, it is sufficient to prove the following proposition.

\begin{prop}\label{DiracforGUnocoeffismallgroupSi}
\begin{equation}
\mu_{H,\pt}: \KK^H_*(\mathcal{S},\pt)\rightarrow\K_H^*(\pt)
\end{equation}
is an isomorphism.
\end{prop}
\begin{proof} It is sufficient to prove
$$
\mathfrak{D}_H=j^H_r(\mathfrak{d}_{H,\mathcal{S}}) \in \KK_0(C^*_r(H,C_{\tau}(\mathcal{S})),C^*_r(H))
$$
is invertible. In fact, we can prove that $\mathfrak{d}_{H, \mathcal{S}}\in \KK^H_0(C_{\tau}(\mathcal{S}),\mathbb{C})$ is invertible. This follows from the fact that $H$ is almost connected amenable together with some formal arguments as follows.

As in the construction in Section \ref{subsection: dual Dirac}, we have the dual Dirac element
$$
\eta_{H,\mathcal{S}}\in \KK^H_0(\mathbb{C},C_{\tau}(\mathcal{S}))
$$
and
\begin{align*}
\mathfrak{d}_{H,\mathcal{S}}\otimes&_{\mathbb{C}}\eta_{H,\mathcal{S}}=1\in \KK^H_0(C_{\tau}(\mathcal{S}),C_{\tau}(\mathcal{S})),\\
 \eta_{H,\mathcal{S}}\otimes&_{C_{\tau}(H)}\mathfrak{d}_{H,\mathcal{S}}=\gamma_{H,\mathcal{S}}\in \KK^H_0(\mathbb{C},\mathbb{C}).
\end{align*}

It is clear that  $H$ is an almost connected amenable group therefore Corollary \ref{Amemnablegammaelement} tells us that
$$
\gamma_{H, \mathcal{S}_H}=1
$$
where $\mathcal{S}_H=H/U\cap H$. By Lemma \ref{restrictionofDiracelements} and Corollary \ref{coro: restriction of gamma} we know that $\gamma_{H, \mathcal{S}_H}=\gamma_{H,\mathcal{S}}$ hence
\begin{equation}
\gamma_{H,\mathcal{S}}=1.
\end{equation}

Now we proved that $\mathfrak{d}_{H,\mathcal{S}}$ hence $\mathfrak{D}_H$ is invertible. As a result we have
$$
\mu_{H,\pt}: \KK^H_*(\mathcal{S},\pt)\rightarrow\K_H^*(\pt)
$$
is an isomorphism. \end{proof}

\begin{proof}[Proof of Proposition \ref{prop: AssembleforGUrealssorbit}] Now it is a corollary of  Proposition \ref{reducetoH} and Proposition \ref{DiracforGUnocoeffismallgroupSi}.\end{proof}

\section{The $G$-orbits on the Flag Variety}\label{section:The G-orbits on the Flag Variety}
We have proved that the assemble map
$$
\mu_{G,\mathcal{O}_{\alpha}^+}: \KK^G_*(S,\mathcal{O}_{\alpha}^+)\rightarrow \K_G^*(\mathcal{O}_{\alpha}^+)
$$ 
is an isomorphism on one $G$-orbit $\mathcal{O}_{\alpha}^+$. In this section we study the geometry of $G$-orbits on $\mathcal{B}$ and in the next section we will piece together orbits.

The result on the $G$-orbits in \cite{mirkovic1992matsuki} is important to our purpose, so we summarize their result here

\begin{thm}\label{morsefunctionorbitsonF}[\cite[1.2, 3.8]{mirkovic1992matsuki}]
Let $G$ be a connected real semisimple Lie group and $\mathcal{B}=G_{\mathbb{C}}/B_{\mathbb{C}}$ as before. Let $U$ be the maximal compact subgroup of $G$.
On the flag variety $\mathcal{B}$ there exists a real value function $f$ such that
\begin{enumerate}
\item $f$ is a Morse-Bott function on $\mathcal{B}$.
\item $f$ is $U$ invariant, hence the gradient flow $\phi: \mathbb{R}\times \mathcal{B}\rightarrow \mathcal{B}$ is also $U$ invariant.
\item The critical point set $\mathcal{C}$ of $f$ consists of finitely many $U$-orbits $\mathcal{O}_\alpha$. The flow preserves the orbits of $G$.
\item The limits $\lim_{t\rightarrow \pm \infty}\phi_t(x):=\pi^{\pm}(x) $ exist for any $x\in \mathcal{B}$. For $\mathcal{O}_\alpha$ a critical $U$-orbit, the stable set
$$
\mathcal{O}_\alpha^+=(\pi^+)^{-1}(\mathcal{O}_\alpha)
$$
is an $G$-orbit, and the unstable set
$$
\mathcal{O}_\alpha^-=(\pi^-)^{-1}(\mathcal{O}_\alpha)
$$
is an $U_\mathbb{C}$-orbit, where $U_\mathbb{C}$ is the complexification of $U$ in $G_\mathbb{C}$.
\item $\mathcal{O}_\alpha^+\cap \mathcal{O}_\alpha^-=\mathcal{O}_\alpha$.
\end{enumerate}
\end{thm}

\begin{coro}
Let $G$ be a connected real semisimple Lie group and $\mathcal{B}=G_{\mathbb{C}}/B_{\mathbb{C}}$ as before. Then the total number of $G$-orbits in $\mathcal{B}$ is finite.
\end{coro}

We will also use the following corollary in \cite{mirkovic1992matsuki}:
\begin{coro}\label{closurerelationoforbits}[\cite[1.4]{mirkovic1992matsuki}]
Let $\mathcal{O}_{\alpha}$ and $\mathcal{O}_{\beta}$ be two critical $U$-orbits. Then the closure $\overline{\mathcal{O}_{\alpha}^+}\supset \mathcal{O}_{\beta}^+$ if and only if
$$
\mathcal{O}_{\alpha}^+\cap \mathcal{O}_{\beta}^-\neq \emptyset.
$$
\end{coro}

From this we can get
\begin{coro}\label{closurerelationofdifferentorbits}
Let $\mathcal{O}_{\alpha}$ and $\mathcal{O}_{\beta}$ be two different critical $U$-orbits, i.e. $\mathcal{O}_{\alpha}\neq \mathcal{O}_{\beta}$. Then $\overline{\mathcal{O}_{\alpha}^+}\supset \mathcal{O}_{\beta}^+$ implies that
the Morse-Bott function $f$ has values
$$
f(\mathcal{O}_{\alpha})>f(\mathcal{O}_{\beta})
$$
\end{coro}
\begin{proof} By Corollary \ref{closurerelationoforbits},
$$
\mathcal{O}_{\alpha}^+\cap \mathcal{O}_{\beta}^-\neq \emptyset.
$$
so there exists an $x\in\mathcal{O}_{\alpha}^+\cap \mathcal{O}_{\beta}^-$.

Since $\lim_{t\rightarrow + \infty}\phi_t(x)\in \mathcal{O}_{\alpha}$, we have
$$
f(\mathcal{O}_{\alpha})\geqslant f(x),
$$
similarly
$$
f(x)\geqslant f(\mathcal{O}_{\beta}).
$$

On the other hand since $\mathcal{O}_{\alpha}$ and $\mathcal{O}_{\beta}$ are critical and  $\mathcal{O}_{\alpha}\neq \mathcal{O}_{\beta}$ we get 
$$
x\not\in \mathcal{O}_{\alpha}, x\not\in \mathcal{O}_{\beta}
$$
hence
$$
f(x)\neq f(\mathcal{O}_{\alpha}), f(x)\neq f(\mathcal{O}_{\beta}).
$$
Therefore we have
$$
f(\mathcal{O}_{\alpha})> f(\mathcal{O}_{\beta}).
$$
\end{proof}

\begin{defi}\label{partialorderonorbitsofF}
We  give a partial order on the set of $G$-orbits of $\mathcal{B}$ which satisfies the following conditions
\begin{enumerate}
\item If $f(\mathcal{O}_{\alpha})>f(\mathcal{O}_{\beta})$, we require that $\mathcal{O}_{\alpha}+>\mathcal{O}_{\beta}+$;

\item If $f(\mathcal{O}_{\alpha})=f(\mathcal{O}_{\beta})$, we choose and fix an arbitrary partial order on them.
\end{enumerate}
\end{defi}

Now let us list all $G$-orbits in $\mathcal{B}$ in ascending order, keep in mind that there are finitely many of them:
\begin{equation}
\mathcal{O}_{\alpha_1}^+<\mathcal{O}_{\alpha_2}^+<\ldots \mathcal{O}_{\alpha_k}^+.
\end{equation}

From the definition we can easily get
\begin{coro}\label{closedsmallorbits}
For any $G$-orbits $\mathcal{O}_{\alpha_i}^+$, the union
$$
\mathcal{Z}_i:=\bigcup_{\mathcal{O}_{\alpha_j}^+\leqslant \mathcal{O}_{\alpha_i}^+}\mathcal{O}_{\alpha_j}^+
$$
is a closed subset of $\mathcal{B}$. Notice that $\mathcal{O}_{\alpha_i}^+\subset \mathcal{Z}_i$
\end{coro}
\begin{proof} It is sufficient to prove that $\mathcal{Z}_i$ contains all its limit points, which is a direct corollary of Definition \ref{partialorderonorbitsofF} and Corollary \ref{closurerelationofdifferentorbits}.
\end{proof}

\begin{rmk}
Corollary \ref{closurerelationofdifferentorbits}, Definition \ref{partialorderonorbitsofF} and Corollary \ref{closedsmallorbits} are not explicitly given in \cite{mirkovic1992matsuki}.
\end{rmk}

\section{The Baum-Connes Conjecture on Flag Varieties}\label{section:Baum-Connes Conjecture on Flag Varieties}
With the construction in Section \ref{section:The G-orbits on the Flag Variety}, we can piece together assemble maps on different orbits.
\begin{prop}\label{sesoforbitsonF}
Let $\mathcal{O}_{\alpha_i}^+$ and $\mathcal{Z}_i$ be as in Theorem \ref {morsefunctionorbitsonF} and Corollary \ref{closedsmallorbits}. Then for $1\leqslant i\leqslant k-1$ we have a short exact sequence of crossed product algebras:
$$
0\rightarrow C^*_r(G, C_0(\mathcal{O}_{\alpha_{i+1}}^+))\rightarrow C^*_r(G, C(\mathcal{Z}_{i+1}))\rightarrow C^*_r(G, C(\mathcal{Z}_{i}))\rightarrow 0.
$$
\end{prop}
\begin{proof} From the construction we also get
\begin{align*}
\mathcal{Z}_i&\subset \mathcal{Z}_{i+1}, ~\mathcal{O}_{\alpha_{i+1}}^+\subset \mathcal{Z}_{i+1},\\
 \mathcal{Z}_i&\cup \mathcal{O}_{\alpha_{i+1}}^+= \mathcal{Z}_{i+1},~ \mathcal{Z}_i\cap \mathcal{O}_{\alpha_{i+1}}^+=\emptyset,
\end{align*}
and $\mathcal{Z}_i$ is closed in $\mathcal{Z}_{i+1}$, $\mathcal{O}_{\alpha_{i+1}}^+$ is open in $\mathcal{Z}_{i+1}$.

Since $\mathcal{B}$ is a compact manifold, we get that  $\mathcal{Z}_i$ and $\mathcal{Z}_{i+1}$ are both compact.

The inclusion gives a short exact sequence of $C^*$-algebras:
\begin{equation}\label{eq: ext without group}
0\rightarrow C_0(\mathcal{O}_{\alpha_{i+1}}^+)\rightarrow C(\mathcal{Z}_{i+1})\rightarrow C(\mathcal{Z}_{i})\rightarrow 0.
\end{equation}

Now we need to go to the reduced crossed-product $C^*$-algebras for which we need the following result.
\begin{lemma}\label{exactcrossedproduct}[\cite[Theorem 6.8]{kirchberg1999permanence}]
Let $G$ be a locally compact group and
$$
0\rightarrow A\rightarrow B\rightarrow C\rightarrow 0
$$
be a short exact sequence of $G$-$C^*$ algebras. Then we have a short exact sequence of reduced cross product $C^*$-algebras:
$$
0\rightarrow C^*_r(G,A)\rightarrow C^*_r(G,B)\rightarrow C^*_r(G,C)\rightarrow 0.
$$
\end{lemma}

Apply Lemma \ref{exactcrossedproduct} to \eqref{eq: ext without group} we get the short exact sequence
\begin{equation}
0\rightarrow C^*_r(G, C_0(\mathcal{O}_{\alpha_{i+1}}^+))\rightarrow C^*_r(G, C(\mathcal{Z}_{i+1}))\rightarrow C^*_r(G, C(\mathcal{Z}_{i}))\rightarrow 0.
\end{equation}
This finishes the proof of Proposition \ref{sesoforbitsonF}.
\end{proof}

From Proposition \ref{sesoforbitsonF} we have the well-known six-term long exact sequence
$$
\begin{CD}
\K^*(C^*_r(G, C_0(\mathcal{O}_{\alpha_{i+1}}^+)))@>>> \K^*(C^*_r(G, C(\mathcal{Z}_{i+1})))@>>> \K^*(C^*_r(G, C(\mathcal{Z}_{i})))\\
@AAA                                     @.                                @VVV \\
\K^{*+1}(C^*_r(G, C(\mathcal{Z}_{i})))@<<< \K^{*+1}(C^*_r(G, C(\mathcal{Z}_{i+1})))@<<< \K^{*+1}(C^*_r(G, C_0(\mathcal{O}_{\alpha_{i+1}}^+))).
\end{CD}
$$
i.e.
\begin{equation}\label{eq1}
\begin{CD}
\K^*_G(\mathcal{O}_{\alpha_{i+1}}^+)    @>>> \K^*_G(\mathcal{Z}_{i+1})@>>>      \K^*_G(\mathcal{Z}_{i})\\
@AAA                             @.                        @VVV \\
\K^{*+1}_G(\mathcal{Z}_{i})@<<< \K^{*+1}_G(\mathcal{Z}_{i+1})@<<< \K^{*+1}_G(\mathcal{O}_{\alpha_{i+1}}^+).
\end{CD}
\end{equation}

Similarly we have
\begin{equation}\label{eq2}
\begin{CD}
\K^{*}_G(C_0(\mathcal{O}_{\alpha_{i+1}}^+)\otimes C_{\tau}(\mathcal{S}))    @>>> \K^*_G(C(\mathcal{Z}_{i+1})\otimes C_{\tau}(\mathcal{S}))@>>>      \K^*_G(C(\mathcal{Z}_{i})\otimes C_{\tau}(\mathcal{S}))\\
@AAA                             @.                        @VVV \\
\K^{*+1}_G(C(\mathcal{Z}_{i})\otimes C_{\tau}(\mathcal{S}))@<<< \K^{*+1}_G(C(\mathcal{Z}_{i+1})\otimes C_{\tau}(\mathcal{S}))@<<< \K^{*+1}_G(C_0(\mathcal{O}_{\alpha_{i+1}}^+)\otimes C_{\tau}(\mathcal{S})).
\end{CD}
\end{equation}

The following proposition claims that \eqref{eq1} and \eqref{eq2} together form a commutative diagram.

\begin{prop}\label{We have six term 3D commuting diagram}
We have the following commutative diagram:
\begin{equation}\label{The six term 3D commuting diagram}
\begin{tikzpicture}[node distance=4.5cm, auto]
  \node (A1) {$\K^{*}_G(C_0(\mathcal{O}_{\alpha_{i+1}}^+)\otimes C_{\tau}(\mathcal{S}))$ };
  \node [right of=A1] (B1) {$\K^*_G(C(\mathcal{Z}_{i+1})\otimes C_{\tau}(\mathcal{S}))$};
  \node [right of=B1] (C1) {$\\K^*_G(C(\mathcal{Z}_{i})\otimes C_{\tau}(\mathcal{S}))$};
  \node [below of=C1, node distance=3cm] (D1) {$\K^*_G(\mathcal{Z}_{i})$};
  \node [left of=D1] (E1) {$\K^*_G(\mathcal{Z}_{i+1})$};
  \node [left of=E1] (F1) {$\K^*_G(\mathcal{O}_{\alpha_{i+1}}^+)$};

  \node (A2) [left of=A1, below of=A1, node distance=2cm] {$\K^{*+1}_G(C_0(\mathcal{O}_{\alpha_{i+1}}^+)\otimes C_{\tau}(\mathcal{S}))$ };
  \node [right of=A2] (B2) {$\K^{*+1}_G(C(\mathcal{Z}_{i+1})\otimes C_{\tau}(\mathcal{S}))$};
  \node [right of=B2] (C2) {$\K^{*+1}_G(C(\mathcal{Z}_{i})\otimes C_{\tau}(\mathcal{S}))$};
  \node [below of=C2, node distance=3cm] (D2) {$\K^{*+1}_G(\mathcal{Z}_{i})$};
  \node [left of=D2] (E2) {$\K^{*+1}_G(\mathcal{Z}_{i+1})$};
  \node [left of=E2] (F2) {$\K^{*+1}_G(\mathcal{O}_{\alpha_{i+1}}^+)$};

\draw[->] (A1) to node {} (B1);
\draw[->] (B1) to node {}(C1);
\draw[->] (C1) to node {}(C2);
\draw[->] (C2) to node {}(B2);
\draw[->] (B2) to node {}(A2);
\draw[->] (A2) to node {}(A1);

\draw[->, densely dotted] (F1) to node {} (E1);
\draw[->, densely dotted] (E1) to node {}(D1);
\draw[->] (D1) to node {}(D2);
\draw[->] (D2) to node {}(E2);
\draw[->] (E2) to node {}(F2);
\draw[->, densely dotted] (F2) to node {}(F1);

\draw[->] (C1) to node {$\mu$}(D1);
\draw[->, densely dotted] (B1) to node {$\mu$}(E1);
\draw[->, densely dotted] (A1) to node {$\mu$}(F1);
\draw[->] (C2) to node {$\mu$}(D2);
\draw[->] (B2) to node {$\mu$}(E2);
\draw[->] (A2) to node {$\mu$}(F2);
\end{tikzpicture}
\end{equation}
where the top and bottom are the six-term exact sequences and the vertical arrows are assemble maps $\mu$.
\end{prop}
\begin{proof} The diagram commutes because all the vertical maps $\mu$ come from the same element
$$
\mathfrak{d}_{G,\mathcal{S}}\in \KK^G_0(C_{\tau}(\mathcal{S}),\mathbb{C})
$$
as in Section \ref{section:The assemble Map on a Single G-Orbit of F}.
\end{proof}

After all these work we are ready to prove Theorem \ref{Baum Connes flag}.

\begin{proof}[Proof of Theorem \ref{Baum Connes flag}] We use induction on the $\mathcal{Z}_i$'s.
First, for $\mathcal{Z}_1=\mathcal{O}_{\alpha_1}^+$, by Proposition \ref{prop: AssembleforGUrealssorbit},
\begin{equation}
\K^*_G(C(\mathcal{Z}_{1})\otimes C_{\tau}(\mathcal{S}))\overset{\mu}{\longrightarrow} \K^*_G (\mathcal{Z}_{1})
\end{equation}
is an isomorphism.

Assume that for $\mathcal{Z}_i$,
\begin{equation}
\K^*_G(C(\mathcal{Z}_{i})\otimes C_{\tau}(\mathcal{S}))\overset{\mu}{\longrightarrow} \K^*_G (\mathcal{Z}_{i})
\end{equation}
is an isomorphism.

By Proposition \ref{prop: AssembleforGUrealssorbit}, the vertical maps on the left face of \eqref{The six term 3D commuting diagram} are isomorphisms. Moreover by induction assumption  the vertical maps on the right face are isomorphism too, hence by a 5-lemma-argument we get the middle vertical maps are also isomorphisms, i.e. for $\mathcal{Z}_{i+1}$,
\begin{equation}
 \K^*_G(C(\mathcal{Z}_{i+1})\otimes C_{\tau}(\mathcal{S}))\overset{\mu}{\longrightarrow} \K^*_G (\mathcal{Z}_{i+1})
\end{equation}
is an isomorphism.

There are finitely many orbits in $\mathcal{B}$ so the induction stops at the largest $\mathcal{Z}_k$ which is $\mathcal{B}$, hence
\begin{equation}
\mu_{G,\mathcal{B}}: \KK^G_*(\mathcal{S},\mathcal{B})\rightarrow \K_G^*(\mathcal{B})
\end{equation}
is an isomorphism.
we finish the proof Theorem \ref{Baum Connes flag}. \end{proof}

\section{An Example: $\SL(2,\mathbb{R})$}\label{section:An Example}
Recall Example \ref{eg: sl2}. If $G=\SL(2,\mathbb{R})$ then $G_{\mathbb{C}}=\SL(2,\mathbb{C})$. We have $\mathcal{S}=G/U=\mathbb{H}$ and
$$
\mathcal{B}=G_{\mathbb{C}}/B_{\mathbb{C}}=\mathbb{C}P^1\cong S^2.
$$
$G_{\mathbb{C}}$ (hence $G$) acts on $\mathcal{B}=\mathbb{C}P^1$ by fractional linear transform
$$
\begin{pmatrix}a &b \\ c & d\end{pmatrix}\cdot \begin{pmatrix}u \\ v\end{pmatrix}:=\begin{pmatrix}au+bv \\ cu+dv\end{pmatrix}.
$$
If we set $z=u/v$, then
\begin{equation}\label{the fractional linear transform}
\begin{pmatrix}a &b \\ c & d\end{pmatrix}\cdot z:= \frac{az+b}{cz+d}.
\end{equation}

From  \eqref{the fractional linear transform} we can see that the action of $G$ on $\mathcal{B}$ is not transitive. In fact, it has three orbits
\begin{align*}
\mathcal{O}_{\alpha_1}^+=& \mathbb{R}\cup \infty\cong S^1 \text{ the equator},\\
\mathcal{O}_{\alpha_2}^+=&\{x+iy|y>0\}\cong \mathbb{C} \text{ the upper hemisphere},\\
\mathcal{O}_{\alpha_3}^+=&\{x+iy|y<0\}\cong \mathbb{C} \text{ the lower hemisphere}.
\end{align*}
$\mathcal{O}_{\alpha_1}^+$ is a closed orbit with dimension $1$; $\mathcal{O}_{\alpha_2}^+$ and $\mathcal{O}_{\alpha_3}^+$ are open orbits with dimension 2.

We look at $\mathcal{O}_{\alpha_1}^+$ first. Take the point $1\in\mathcal{O}_{\alpha_1}^+$. The isotropy group at $1$ is the upper triangular group $B$ in $\SL(2,\mathbb{R})$. So
$$
\K^*_G(\mathcal{O}_{\alpha_1}^+)=\K^*_B(pt).
$$

$B$ is solvable hence amenable and $\mathbb{Z}/2\mathbb{Z}$ is the maximal compact group of $B$. By Theorem \ref{Diracforamenablegroups}
$$
\K^0_B(pt)=R(\mathbb{Z}/2\mathbb{Z})
$$
is the representation ring of the group with two elements and
$$
\K^1_B(pt)=0.
$$
So
$$
\K^0_G(\mathcal{O}_{\alpha_1}^+)=R(\mathbb{Z}/2\mathbb{Z})
$$
and
$$
\K^1_G(\mathcal{O}_{\alpha_1}^+)=0.
$$

For $\mathcal{O}_{\alpha_2}^+$ and $\mathcal{O}_{\alpha_3}^+$, the isotropy groups are both
$$
T=\left\{\begin{pmatrix}\cos\theta &\sin\theta \\ -\sin\theta & \cos\theta\end{pmatrix}\right\}
$$
hence by the similar reason to  $\mathcal{O}_{\alpha_1}^+$ we have
$$
\K^0_G(\mathcal{O}_{\alpha_2}^+)=\K^0_G(\mathcal{O}_{\alpha_3}^+)=\K^0_T(pt)=R(T)
$$
is the representation ring of $T$ and
$$
\K^1_G(\mathcal{O}_{\alpha_2}^+)=\K^1_G(\mathcal{O}_{\alpha_3}^+)=\K^1_T(pt)=0.
$$

Now $\mathcal{O}_{\alpha_2}^+\cup\mathcal{O}_{\alpha_3}^+$ is open in $\mathcal{B}$ so as in the last section we have the short exact sequence
\begin{equation}
0\longrightarrow C_0(\mathcal{O}_{\alpha_2}^+\cup\mathcal{O}_{\alpha_3}^+)\longrightarrow C(\mathcal{B})\longrightarrow C(\mathcal{O}_{\alpha_1}^+)\longrightarrow 0
\end{equation}
and further
$$
0\longrightarrow C_r^*(G,\mathcal{O}_{\alpha_2}^+\cup\mathcal{O}_{\alpha_3}^+)\longrightarrow C_r^*(G,\mathcal{B})\longrightarrow C_r^*(G,\mathcal{O}_{\alpha_1}^+)\longrightarrow 0.
$$
i.e.
\begin{equation}
0\longrightarrow C_r^*(G,\mathcal{O}_{\alpha_2}^+)\oplus C_r^*(G,\mathcal{O}_{\alpha_3}^+)\longrightarrow C_r^*(G,\mathcal{B})\longrightarrow C_r^*(G,\mathcal{O}_{\alpha_1}^+)\longrightarrow 0.
\end{equation}
We get the six-term exact sequence
\begin{equation}
\begin{CD}
\K^0_G(\mathcal{O}_{\alpha_2}^+) \oplus\K^0_G(\mathcal{O}_{\alpha_3}^+)    @>>> \K^0_G(\mathcal{B})@>>>      \K^0_G(\mathcal{O}_{\alpha_1}^+)\\
@AAA                                                  @.                        @VVV \\
\K^1_G(\mathcal{O}_{\alpha_1}^+)@<<< \K^1_G(\mathcal{B})@<<< \K^1_G(\mathcal{O}_{\alpha_2}^+) \oplus\K^1_G(\mathcal{O}_{\alpha_3}^+) .
\end{CD}
\end{equation}
Combine with the previous calculation we get
\begin{equation}
0\longrightarrow R(T)\oplus R(T)\longrightarrow \K^0_G(\mathcal{B})\longrightarrow R(\mathbb{Z}/2\mathbb{Z})\longrightarrow 0
\end{equation}
and $\K^1_G(\mathcal{B})=0$.

In conclusion we have
\begin{equation}\label{the K-theory of SL(2,R) on F}
\begin{split}
\K^0_G(\mathcal{B})\approx& R(T)\oplus R(T)\oplus R(\mathbb{Z}/2\mathbb{Z}),\\
\K^1_G(\mathcal{B})=& 0.
\end{split}
\end{equation}

Next we look at $\KK^G_*(\mathcal{S},\mathcal{B})$. By Corollary \ref{coro: reduce to max compact subgroup} we have $\KK^G_*(\mathcal{S},\mathcal{B})\cong \K^*_U(\mathcal{B})$. We know that for $G=\SL(2,\mathbb{R})$ the maximal compact subgroup $U=T$.
By Bott periodicity we have
\begin{equation}
\K^0_U(\mathcal{O}_{\alpha_2}^+)=\K^0_U(\mathcal{O}_{\alpha_3}^+)\cong \K^0_U(\mathbb{C})\cong \K^0_U(pt)= R(U)=R(T)
\end{equation}
and
\begin{equation}
\K^1_U(\mathcal{O}_{\alpha_2}^+)=\K^1_U(\mathcal{O}_{\alpha_3}^+)\cong \K^1_U(\mathbb{C})\cong \K^1_U(pt)=0.
\end{equation}

As for $\mathcal{O}_{\alpha_1}^+$, we notice that $U$ acts on $\mathcal{O}_{\alpha_1}^+\cong S^1$ by "square", so the isotropy group is $\mathbb{Z}/2\mathbb{Z}$. Hence
$$
\K^0_U(\mathcal{O}_{\alpha_1}^+)=R(\mathbb{Z}/2\mathbb{Z}) \text{ and } \K^1_U(\mathcal{O}_{\alpha_1}^+)=0.
$$
By the  six-term long exact sequence we have
\begin{equation}\label{the K-theory of U=T on F}
\begin{split}
\K^0_U(\mathcal{B})\approx& R(T)\oplus R(T)\oplus R(\mathbb{Z}/2\mathbb{Z}),\\
\K^1_U(\mathcal{B})=& 0.
\end{split}
\end{equation}

\eqref{the K-theory of SL(2,R) on F} and \eqref{the K-theory of U=T on F} is compatible with the 
 Baum-Connes conjecture (in fact, Connes-Kasparov conjecture as in Section \ref{subsection: spin case}) which states that
\begin{equation}
\K^*_U(\mathcal{B})\cong \K^*_G(\mathcal{B}).
\end{equation}

\begin{rmk}
Using Bott periodicity theorem we can obtain precisely the algebra structure of $\K^*_U(\mathcal{B})$ as in \cite{segal1968equivariant}. Therefore Baum-Connes conjecture will be a powerful tool to investigate $\K_G(\mathcal{B})$ and to study the representation theory of $G$.
\end{rmk}

\bibliography{BCfvbib}{}
\bibliographystyle{plain}
\end{document}